\documentclass[11pt]{article}

\input{amssym.def}
\input{amssym}
\usepackage{eufrak}

\newtheorem{theorem}{Theorem}
\newtheorem{proposition}[theorem]{Proposition}

\def\qq{q^{-1}}

\def\T{{\cal T}}
\def\L{{\cal L}}

\def\ot{\otimes}
\def\C{{\Bbb C}}
\def\End{{\rm End}}
\def\vv{V^{\otimes 2}}
\def\RR{{\cal R}}

\def\al{{\alpha}}

\def\be{\begin{equation}}
\def\ee{\end{equation}}

\begin{document}

\makeatletter
\renewcommand{\theequation}{{\thesection}.{\arabic{equation}}}
\@addtoreset{equation}{section} \makeatother

\title{Centers in Generalized  Reflection Equation algebras}

\author{\rule{0pt}{7mm} Dimitri
Gurevich\thanks{gurevich@ihes.fr}\\
{\small\it LAMAV, Universit\'e de Valenciennes,
59313 Valenciennes, France}\\
{\small \it and}\\
{\small \it
Interdisciplinary Scientific Center J.-V.Poncelet}\\
{\small \it Moscow 119002, Russian Federation}\\
\rule{0pt}{7mm} Pavel Saponov\thanks{Pavel.Saponov@ihep.ru}\\
{\small\it
National Research University Higher School of Economics,}\\
{\small\it 20 Myasnitskaya Ulitsa, Moscow 101000, Russian Federation}\\
{\small \it and}\\
{\small \it
Institute for High Energy Physics, NRC "Kurchatov Institute"}\\
{\small \it Protvino 142281, Russian Federation}}

\maketitle

\begin{abstract}
In the Reflection Equation (RE) algebra associated with an involutive or Hecke symmetry $R$ the center is generated by elements ${\rm Tr}_R L^k$
(called the quantum power sums), where $L$ is the generating matrix of this algebra and  ${\rm Tr}_R$ is the $R$-trace associated with $R$. We
consider the problem: whether it is so in certain RE-like algebras depending on spectral parameters. Mainly, we deal with algebras similar to those considered
in \cite{RS} (we call them algebras of RS type). These algebras are defined by means of some current (i.e. depending on parameters) $R$-matrices arising from involutive and Hecke
symmetries via the so-called Baxterization procedure. We define quantum  power sums in the algebras of RS  type and show that  the lowest quantum
power sum in such an algebra is central iff the value of the "charge" $c$ entering its definition  is critical. We exhibit  the dependance of this critical value on the bi-rank of the
initial symmetry $R$. Besides, we show that if the bi-rank of $R$ is $(m|m)$, and the value of $c$ is critical, then all  quantum
power sums are central.
\end{abstract}

{\bf AMS Mathematics Subject Classification, 2010:} 81R50

{\bf Keywords:} Reflection equation algebras, algebras of Reshetikhin-Semenov-Tian-Shansky type, charge, quantum powers of the generating matrix,  quantum power sums

\section{Introduction}

The best known Quantum Matrix (QM) algebras related to Quantum Groups (QG) $U_q(sl(N))$ are the so-called RTT and Reflection Equation (RE) ones.
In general, any  QM algebra  is defined via the {\em generating matrix} with entries belonging to this algebra and subject to a system involving
an $R$-matrix or a couple of them (see \cite{IOP}).

By definition, an $R$-matrix is an operator $R:\vv\to \vv$ subject to  the so-called {\em braid relation} (also called Quantum Yang-Baxter Equation)
\be
R_{12}\, R_{23}\, R_{12}= R_{23}\,R_{12}\, R_{23}, \qquad R_{12}=R\ot I,\quad  R_{23}=I\ot R.
\label{braid}
\ee
Hereafter, $V$ is a vector space  of dimension $N$ over the ground field $\C$, and $I$ is the identity operator. Also, as usual, the low indices indicate the
position where a matrix or an operator is located. Below, such operators $R$ will be also referred to as {\em braidings}.

The QM algebras associated with  braidings will be called {\em constant} by contrast with {\em current} algebras considered below: the defining relations of the
latter algebras contain some functions depending on spectral parameters.

Recall that the constant RTT and RE algebras are respectively defined by the following
systems of relations on their generators $t_i^j$ and $l_i^j$:
$$
R\, T_1\, T_2- T_1\, T_2\,R=0,\quad  T=\|t_i^j\|_{1\leq i, j \leq N},
$$
\be
R\, L_1\, R \, L_1- L_1\,R\, L_1\,R=0,\quad  L=\|l_i^j\|_{1\leq i, j \leq N}.
\label{RE}
\ee
We denote these algebras $\T(R)$ and $\L(R)$ respectively.

All braidings $R$  are assumed  to be {\em skew-invertible} (see section 2) involutive or Hecke symmetries. An involutive symmetry obeys an additional
relation $R^2=I$, while a Hecke symmetry is subject to the relation
$$
(R-q I)(R+\qq I)=0,\quad q\in \C\setminus \{0, \pm 1\},
$$
where $q$ is assumed to be {\it generic}, that is $q^k\not= 1$ for any integer $k$. As a consequence, $k_q=\frac{q^k-q^{-k}}{q-\qq}\not=0$.

Note that in the both QM algebras  analogs of basic {\em symmetric polynomials} (elementary, full, Schur polynomials, power sums) are well-defined.
However, the properties of these quantum symmetric polynomials  in the mentioned QM algebras differ drastically. In the RTT algebras the quantum symmetric polynomials
generate {\it a commutative subalgebra}. In the RE algebras they are  central.

We are especially interested in quantum analogs of the {\em  power sums}. Note that in the algebras $\L(R)$ the {\em quantum
power sums}\footnote{This term is motivated by the following reason: if a matrix $L$ has commutative entries and $\mu_i$ are its eigenvalues, then
${\rm Tr}\, L^k= \sum_k \mu_i^k$.} look like their classical counterparts, namely, ${\rm Tr}_{R}L^k$, $k\ge 1$, whereas  in the algebras $\T(R)$ they cannot
be presented in such a simple form. Hereafter, the notation  ${\rm Tr}_R$ stands for  the so-called $R$-trace, which is associated
with any skew-invertible braiding $R$ (see section 2).

 Note that the term
``RE algebra" appeared in connection with models describing dynamics of particles in domains with boundaries (see \cite{C, S, KS}). A method of constructing
constant RE algebras via couples of RTT algebras was suggested in \cite{FRT}. A generalization of this construction to some current algebras gave rise to the
RE-like algebras  from \cite{RS} (section 2, formula (11)).

The main objective of the present note is the following: we study in what current RE-like algebras the centrality  of the quantum power sums is also valid.
Namely,  we consider RE-like algebras of general form. Each of them is defined by the following system
\be
(R+g_1(u,v)I) L_1(u)(R+g_2(u,v)I)L_1(v)=L_1(v)(R+g_3(u,v)I) L_1(u) (R+g_4(u,v)I),
\label{ge}
\ee
where $R$ is a skew-invertible involutive or Hecke symmetry  and $g_k(u,v),\, k=1,2,3,4$  are some functions. We call the algebras defined by  (\ref{ge}) {\em the generalized
reflection equation algebras}.

Note that the factors $R+g_k(u,v)I$ are usually assumed to be current $R$-matrices (may be, with a charge introduced inside). In this paper we mainly deal with
algebras similar to those from \cite{RS}, but associated with other current $R$-matrices, which arise from involutive and Hecke symmetries via the Baxterization
procedure. We call them {\em the algebras of RS type}.

Given  an algebra of RS type,  we define {\em quantum matrix powers} of its {\em generating matrix} $L(u)$ via the following formula
\be
L^{[k]}(u)=L_1(c^{k-1}u)\,L_1(c^{k-2}u)\dots L_1(cu)\,L_1(u),\,\,\, k=1,2,...
\label{string}
\ee
The factor $c$ here is called the {\em charge}.

Observe that the notion of quantum powers of  matrices was introduced by Talalaev (see \cite{T}) in connection with the rational
Gaudin model. In another way  quantum matrix powers can be defined in  the Drinfeld's Yangians and their generalizations (see \cite{GS, GST}).
In all mentioned algebras this notion is motivated by the corresponding versions of the Cayley-Hamilton identity. Our present version
of quantum matrix powers in the algebras of RS type  is motivated by the commutation properties of these quantum powers (see section 4).

Observe that if $R$ is an involutive symmetry, the quantum power $L^{[k]}(u)$ is defined by an analogous formula but the shifts of the argument $u$ are additive:
$u\mapsto u+k\,c$, $k\ge1$.

By using  the quantum  power sums we define quantum analogs of power sums  as follows ${\rm Tr}_R L^{[k]}(u)$.

Our first main result is the following.
We show that in any algebra of RS type the first quantum power sum ${\rm Tr}_R L(u)$ is central, iff the value of the charge $c$ is {\em critical}. This result is similar to a claim
from \cite{RS}, where it is established for the algebras related to the QG $U_q(g)$ ($g$ is a simple classical Lie algebra).  However, if in the  setting  from \cite{RS} the critical value of the charge is related to the dual Coxeter number of $g$,  in our setting this value is determined by the bi-rank of the initial symmetry. In the case $g=sl(N)$ (below,
it will be called  {\em standard})   our result is a generalization of that from \cite{RS}.

Stress that the method used in the present note differs drastically  from that of \cite{RS}.  Namely, using  the so-called $R$-matrix technique\footnote{By
$R$-matrix technique we mean a collection of formulae valid for all (skew-invertible)  braidings or symmetries regardless their concrete forms (see section 2). Note that this technique does not use any object of QG type.} we succeeded to generalize the mentioned result from  \cite{RS} to algebras associated with symmetries (involutive or Hecke) of general type. It should be emphasized that the family of such symmetries is large  enough. A way of  constructing such symmetries for all possible bi-ranks was exhibited in \cite{G}. (see footnote
\ref{fooo}).

Besides, we show, that  if the
bi-rank of  the initial symmetry $R$ is $(m|m)$ then all quantum power sums are central in the corresponding algebras of RS type for the critical value of 
the charge\footnote{Note that in this case the critical value of
the charge is 1.}. Unfortunately, we have not succeeded in establishing (or rejecting)  this property if the bi-rank of $R$ is arbitrary. Nevertheless, we want to stress
usefulness and powerfulness of the $R$-matrix technique for such sort of study.

By concluding Introduction we want to observe
 that we do not expand the generating matrices of the algebras under consideration in Laurent series, as usually done in the (double) Yangians.
Thus, we avoid the problem of a normal ordering of quantum powers of matrices.  Note that in general such an ordering is not well-defined and can be only done
in algebras arising from a  couple of RTT algebras in the spirit of (\ref{pm}). We refer the reader to the paper \cite{KM}, where such a construction in a quantum affine vertex algebra
associated with trigonometrical $R$-matrix, is exhibited, but no RE-like algebra is considered.

\medskip

\noindent
{\bf Acknowledgement}  The work of P.S. has been funded by the Russian Academic Excellence Project `5-100' and  also partially supported by the RFBR grant 16-01-00562.

\section{Symmetries and QM algebras}

In this section we consider a relation between constant RTT and RE algebras associated with symmetries. However, first, we recall some elements of the so-called
$R$-matrix technique. For detail the reader is referred to \cite{GPS, O}.

Let $R:\vv\to \vv$ be  a Hecke symmetry. Let us consider the associated $R$-skew-symmetric algebra
\be
\mbox{\Large $\mathsf{\Lambda}$}_R(V)=T(V)/\langle\mathrm{ Im}(\qq I+R)\rangle,
\label{skew}
\ee
where $T(V)$ is the free tensor algebra of the space $V$ and $\langle I\rangle$ stands for the ideal generated by a subset $I$. Let $\mbox{\large $\mathsf{\Lambda}$}_R^{(k)}(V)$
be its $k$-th order homogenous component and
$$
P_{{\mathsf{\Lambda}}_R}(t)=\sum_{k\geq 0} t^k \,\mathrm{dim}\, \mbox{\large $\mathsf{\Lambda}$}_R^{(k)}(V)
$$
be the corresponding  Poincar\'e-Hilbert series. If $R$ is involutive, then we put $q=1$ in formula (\ref{skew}).

It is known (see \cite{GPS}) that for any symmetry $R$ (involutive or Hecke)  the series $P_{{\mathsf{\Lambda}}_R}(t)$ is a rational function.
We say that the bi-rank of $R$ is $(m|n)$ if $P_{\mathsf{\Lambda}_R}(t)$ is a ratio of two coprime polynomials and $n$ and $m$ are degrees of the numerator
and of the denominator respectively.

Also, we assume $R$ to be skew-invertible. This means that  there exists an operator $\Psi:\vv\to \vv$ such that
$$
\mathrm{ Tr}_{(2)} R_{12}\Psi_{23}=P_{13}.
$$
Hereafter, $P$ stands for the usual flip or its matrix. Define the following operators
$$
B=\mathrm{Tr}_{(1)} \Psi,\quad C=\mathrm{Tr}_{(2)} \Psi
$$
and introduce the {\it $R$-trace} by the rule
$$\mathrm{Tr}_R A=\mathrm{Tr} (C\cdot A).
$$
Here $\mathrm{Tr}$ is the usual trace and $A:V \to V$ is an $N\times N$ matrix. Note that we fix a basis in $V$ and consequently identify the operators and matrices.

As  for the operator $B$, it comes in the definition of a product in the algebra $\End(V)$, in a sense coordinated with the initial symmetry $R$. The operators $B$ and $C$ are invertible
and connected  by the following formula
\be
B\cdot C=q^{2(n-m)}\,I,
\label{BC}
\ee
provided $R$ is a skew-invertible  Hecke symmetry of bi-rank $(m|n)$. Also, note that for the $R$-trace the following formulae hold true
\be
\mathrm{ Tr}_{R(2)} R_{12}=I_1, \qquad \mathrm{ Tr}_{R} I=q^{n-m} (m-n)_q.
\label{ff}
\ee

Also, we need the following  formula proved in \cite{O}
\be
\mathrm{ Tr}_{R(2)}R_{12} A_1 R^{-1}_{12}=\mathrm{ Tr}_{R(2)}R_{12}^{-1} A_1 R_{12} =I_1 (\mathrm{Tr}_{R} A),
\label{Oleg}
\ee
where $A$ is any $N\times N$ matrix (may be with non-commutative entries).

With the use of (\ref{Oleg}) one can easily show that the quantum power sums $\mathrm{Tr}_R L^k$ belong to the center of the algebra $\L(R)$.
Indeed, upon applying  the defining relations of the algebra $\L(R)$ $k$  times, we arrive to the following equality
$$
R_{12} L_1 R_{12} (L_1)^k= (L_1)^k R_{12} L_1R_{12},\qquad k\ge 1.
$$
Let us multiply this relation by $R^{-1}$ from the both sides and apply $\mathrm{Tr}_R$ at the second position. Then we get
$$
\mathrm{Tr}_{R(2)} L_1 R_{12}(L_1)^k R_{12}^{-1}= \mathrm{Tr}_{R(2)} R_{12}^{-1}(L_1)^k R_{12} L_1.
$$

Now, by taking in consideration (\ref{Oleg}), we arrive to the final conclusion:
$$
L\, (\mathrm{Tr}_R L^k)=(\mathrm{Tr}_R L^k)\, L,
$$
meaning that $\mathrm{Tr}_R L^k$ belong to the center of the algebra  $\L(R)$ for any $k$.

It should be emphasized that this method does not work in the algebras $\T(R)$. Though quantum analogs of power sums can also be defined in these algebras, they cannot be
reduced to the form $\mathrm{Tr}_{R} T^k$ (the explicit form of the power sums in this case can be found for instance in \cite{IOP}). Besides, these elements are
not central. However, they generate a commutative subalgebra of $\T(R)$. Thus, we can conclude that algebraic properties of the QM algebras $\L(R)$ and $\T(R)$
differ drastically. Besides, their coalgebraic structures differ as well, but we do not consider them here.

Now, let us recall the method from  \cite{FRT} enabling one to construct an RE algebra via a couple of the RTT ones. However, we reproduce it in a much more general setting by assuming
$R$ to be  any skew-invertible symmetry.

Consider an algebra, composed of two copies of the RTT algebras, $T^+$ and $T^-$ being their generating matrices
\be
R\, T^{\pm}_1 T^{\pm}_2=T^{\pm}_1 T^{\pm}_2 R,\quad R\, T^{+}_1 T^{-}_2=T^{-}_1T^{+}_2 R.
\label{pm}
\ee
The latter equality defines the {\em permutation relations } among the generators of the RTT copies.

Now, assume  the symmetry $R$  to be of bi-rank $(m|0)$\footnote{\label{fooo}As an example of such a Hecke symmetry we mention that coming from the QG $U_q(sl(N))$. In
this case $m=N$. However, in general it is not so. Thus, in \cite{G} for any $N\geq 2$ there are constructed involutive and Hecke symmetries of bi-rank $(2|0)$ acting
in the space $\vv$, where  $\dim V=N$. For the corresponding $R$-skew-symmetric algebra the Poincar\'e-Hilbert series is $P_{\mathsf{\Lambda}_R}=1+Nt+t^2$.
Also, in \cite{G} there is exhibited a method of "gluing" of Hecke symmetries. Thus, by gluing a Hecke symmetry of bi-rank $(m|0)$ and that of  bi-rank $(0|n)$, we get a Hecke symmetry of
bi-rank $(m|n)$.}.
Then an $R$-determinant $\mathrm{det}_R T$ is well-defined in the RTT algebra $\T(R)$ (if $R=P$ is the usual flip, $\mathrm{det}_R T$ turns into the usual
determinant). For defining an antipode we assume $\mathrm{det}_R \, T$ to be central. Then the quotient-algebra $\T(R)/\langle \mathrm{det}_RT-1 \rangle$  is
well-defined and an antipode $S$ in this algebra is also well-defined\footnote{If the $R$-determinant is not central, it is possible to deal with the localization
$(\mathrm{det}_R T)^{-1}\T(R)$. In this algebra an antipode exists too. An example of a symmetry $R$ such that the determinant $det_R \, T$  is not central is
exhibited in \cite{GS}.}. By using the relation $T\,S(T)=S(T)\, T=I$ we can rewrite the last of the relations  (\ref{pm}) as follows
$$
S(T^-_1) R_{12} T^+_1=T^+_2 R_{12} S(T_2^-).
$$
Now, it is  easy  to prove the following proposition.

\begin{proposition}
The matrix $L=T^+S(T^-)$ meets the system {\rm (\ref{RE})}.
\end{proposition}

In the standard case this claim was established in \cite{FRT}. Observe that in this case one usually imposes some complimentary conditions on the matrices $T^{\pm}$ by
assuming them to be, respectively, upper and lower  triangular and imposing additional relations on the diagonal elements. These conditions enable one to reduce the algebra
(\ref{pm}) to the "classical size".

Let us mention two properties of the  algebra $\L(R)$ while $R$ is standard. First, this algebra is a deformation of the commutative algebra $\mathrm{Sym}(gl(N))$, i.e. as $q\to 1$ the
 algebra $\L(R)$ turns into the latter one and for a generic $q$ the dimensions of homogeneous components of both algebras are equal to each other. Observe that the same property
is also valid for the corresponding RTT algebra $\T(R)$. Second, the algebra $\L(R)$ is covariant with respect
to the adjoint action of the QG  $U_q(gl(N))$. Whereas, for the corresponding RTT algebra $\T(R)$ it is not so.

\section{Centrality of $\mathrm{Tr}_R L(u)$}

Let us consider algebras generated by  matrices $L(u)$ subject to the  system (\ref{ge}), where $R$ is assumed to be  a skew-invertible involutive or Hecke symmetry.
Note that if all $g_k$ are trivial, we recover the constant RE algebras. Now, let us suppose that none of $g_k$ is trivial. However, for the moment we do not specialize
the functions $g_k=g_k(u,v),\, k=1,2,3,4$.

Below, we employ the method which has been used in section 2 for studying the center of the  algebra $\L(R)$. Namely, we
apply the $R$-trace $\mathrm{Tr}_{R(2)}$ to the both sides of (\ref{ge}). Then in the left hand side we have
$$ \mathrm{Tr}_{R(2)}(R\, L_1(u)\, ((q-\qq)\, I+R^{-1}) \, L_1(v)+g_1\,L_1(u)\, (\mathrm{Tr}_{R(2)}\,R)\,L_1(v)+g_2\,(\mathrm{Tr}_{R(2)}\, R)\,L_1(u)\,L_1(v)+$$
\be (\mathrm{Tr}_{R(2)}\, I_2)\,g_1\,g_2\,L_1(u)\,L_1(v))
=\mathrm{Tr}_R(L_1(u))L_1(v)+  \left((q-\qq)+g_1+g_2+\al g_1\,g_2\right) L_1(u)\, L_1(v), \label{last} \ee
where $\al=q^{n-m} (m-n)_q$.  Here, we used  formulae (\ref{Oleg}), (\ref{ff}) and  that   $R=(q-\qq)I+R^{-1}$.
Hereafter, we are mainly dealing with Hecke symmetries $R$. 

If $R$ is a  skew-invertible involutive symmetry,  this formula is still valid but the term $q-\qq$ vanishes and
$\al=m-n$.

By applying the same method in the right hand side of (\ref{last}), we arrive to the expression
$$L_1(v) \mathrm{Tr}_R(L_1(u))+ \left((q-\qq)+g_3+g_4+\al g_3\,g_4\right) L_1(v)\, L_1(u).$$
Hereafter, $\mathrm{Tr}_{R}(L_1(u))=\mathrm{Tr}_{R(1)}(L_1(u))$.
It is evident that  the relation
\be \mathrm{Tr}_R(L_1(u))L_1(v)=L_1(v) \mathrm{Tr}_R(L_1(u)), \label{cent}\ee
meaning that the lowest quantum power sum $\mathrm{Tr}_R\,L_1(u)$ is central, is valid if
the following holds
\be (q-\qq)+g_1+g_2+\al\, g_1\, g_2=0,\,\,(q-\qq)+g_3+g_4+\al\, g_3\, g_4=0. \label{sums} \ee
Note that this condition is also necessary for (\ref{cent}) since any nontrivial linear combination of the products $L_1(v)\, L_1(u)$ and $L_1(u)\, L_1(v)$ differ from 0 for generic $u$ and $v$.

Usually, one supposes that the factors $R+g_k(u,v)$ entering the relation (\ref{ge}) are current  $R$-matrices
may be with shifted parameters. The current $R$-matrices can be constructed via the Baxterization procedure. The following claims are known. Note that for the first time such a Baxterization procedure of Hecke symmetries of general type was performed in \cite{G}.

\begin{proposition} 1. Let $R$ be an involutive symmetry. Then the sum\rm
\be
R(u,v) = R+\frac{a}{u-v}\, I,\quad a\in \C
\label{inv}
\ee
\it
 is a current $R$-matrix, i.e. it meets the Quantum Yang-Baxter Equation with spectral parameters
$$
R_{12}(u,v) \, R_{23}(u,w)\, R_{12}(v,w) = R_{23}(v,w) R_{12}(u,w) \, R_{23}(u,v).
$$

2. Let $R=R_q$ be a Hecke symmetry. Then\rm
\be
R(u,v) = R_q-\frac{(q-\qq) u}{u-v}\, I
\label{Hec}
\ee\it
 is a current $R$-matrix.
\end{proposition}

Below, without loss of generality we set $a=1$ in (\ref{inv}).

The current $R$-matrices (\ref{inv}) and (\ref{Hec}) are called respectively rational and trigonometrical.
Upon replacing $u\to q^{-{2u}}$, $v\to q^{-{2v}}$ in (\ref{Hec}), we get the
following form of the trigonometrical $R$-matrix
\be
R_q+\frac{q-\qq}{q^{2(u-v)}-1}\, I=R_q+\frac{q^{-(u-v)}}{({u-v})_q}\, I.
\label{tri}
\ee
It tends to the rational $R$-matrix as $q\to 1$.

Now, let us write a trigonometrical $R$-matrix in the form
\be
R(u,v)=\RR(u/v),\quad \RR(x)=R+f(x)\, I,
\label{RRR}
\ee
where $f(x)$ is defined by
\be
f(x)=-\frac{(q-\qq)x}{x-1}.
\label{fff}
\ee

Consider a particular case of the system (\ref{ge}) corresponding to the following choice
\be
g_1=f(u/v),\quad g_2=f(vc/u),\quad g_3=f(uc/v),\quad g_4=f(v/u),
\label{ggg}
\ee
where $f(x)$ is defined by (\ref{fff}). The number $c\in\C$ is called the {\em charge}. Then the defining relations (\ref{ge}) can be presented as
\be
\RR(u/v)\, L_1(u)\, \RR(vc/u)\, L_1(v)=L_1(v)\,\RR(uc/v)\,L_1(u)\,\RR(v/u).
\label{RS}
\ee

In the standard case this system is equivalent to  that from \cite{RS}. We have only to note that our charge is $c=e^{hC}$, where $C$ is the charge used in
\cite{RS}. Besides, our normalization of the $R$-matrix differs from that from \cite{RS} by the factor $P$ (the usual flip).

\begin{proposition}
In the algebra introduced by {\rm (\ref{RS})} with the $R$-matrix $R(u,v)$ defined by {\rm (\ref{RRR})} and {\rm (\ref{fff})} the first quantum power sum
$\mathrm{Tr}_R L(u)$ is central iff the value of the charge $c$ equals $q^{2(m-n)} $. 
\end{proposition}

This value of $c$ is called {\em critical}.

{\bf Proof.} It suffices to check that the conditions (\ref{sums}) with the functions $g_k$ defined by (\ref{ggg}) and the function $f(x)$ defined by (\ref{fff}) is valid iff the value of
$c$ is critical. \hfill \rule{6.5pt}{6.5pt}

\medskip

Thus, we can conclude  that the family of the generalized RE algebras with central element $\mathrm{Tr}_R  L(u)$ is sufficiently large. For instance, as follows from the conditions (\ref{sums}),
the centrality of this element remains valid if we interchange  either the functions $g_1$ and $g_2$ or $g_3$ and $g_4$, or if we simultaneously apply the both interchanges
in the defining relations (\ref{ge}).

Any algebra defined by (\ref{RS}) with an arbitrary skew-invertible Hecke symmetry $R$ and with $f(x)$ defined by (\ref{fff}) will be called {\em an algebra of RS type}.

If the initial symmetry $R$ is involutive, we set
\be
\RR(x)=R+f(x)\, I,\quad  f(x)=\frac{1}{x}.
\label{ra}
\ee
Then, the rational version of the RS type algebras is defined by the following system
\be
\RR(u-v)\, L_1(u)\, \RR(v-u+c)\, L_1(v)=L_1(v)\, \RR(u-v+c)\,L_1(u)\, \RR(v-u).
\label{ratRS}
\ee
The summand $c$ is also called the {\em charge}. In this case the element $\mathrm{Tr}_R L(u)$ is central if and only if $c=n-m$. This is the {\em critical
value} of the charge in the rational case.

\section{Study of higher quantum power sums}

In this section we deals with the algebras of RS type defined by the system (\ref{RS}). 

Let us consider the quantum matrix power  $L^{[k]}(v)$ defined by   (\ref{string}). We state that in the product
$$
L^{[k]}_1(v)\,\RR(uc/v) L_1(u)\,\RR(v/u)
$$
the quantum matrix power $L^{[k]}(v)$ can be pushed forward to the outermost right position modulo the defining relations in this algebra. Namely,  this property  motivates our definition
of quantum matrix powers.

Let us verify this claim for the case $k=2$. We have
\begin{eqnarray*}
L_1(cv)L_1(v)\RR(uc/v) L_1(u)\RR(v/u)\!\!\!&=&\!\!\! L_1(cv)\RR(u/v) L_1(u)\RR(vc/u) L_1(v)\\
&=&\!\!\! \RR(u/cv)L_1(u)\RR(c^2 v/u) L_1(cv)L_1(v).
\end{eqnarray*}

The second equality is valid in virtue of the relation (\ref{RS}), where $v$ is replaced by $cv$. By pursuing this procedure
 we arrive to a series of  relations
$$
L^{[k]}_1(v)\,\RR(uc/v) L_1(u)\,\RR(v/u)=\RR(u/c^{k-1}v)\,L_1(u)\,\RR(c^k v/u)\,L^{[k]}_1(v),\quad  k\ge 2.
$$

We multiply this equality  by $\RR^{-1}(u/c^{k-1}v)$ from the left and by  $\RR^{-1}(v/u)$ from the right. Then  we get
$$
\RR^{-1}(u/c^{k-1}v)\, L^{[k]}_1(v)\,\RR(uc/v)\, L_1(u)=L_1(u)\,\RR(c^k v/u)\,L^{[k]}_1(v)\, \RR^{-1}(v/u).
$$
Now, taking into account that
$$
(R+g\,I)^{-1}=\frac{R-g\,I-(q-\qq)\,I}{1-g(g+q-\qq)}=\frac{(R^{-1}-g\,I)}{1-g(g+q-\qq)},
$$
for any function $g$, we arrive to the following equality
$$
\frac{(R^{-1}-f(u/c^{k-1}v)\,I)\,L^{[k]}_1(v)(R+f(uc/v)\,I)\, L_1(u)}{1-f(u/c^{k-1}v)(f(u/c^{k-1}v)+q-\qq)}=$$
$$\frac{ L_1(u)\,(R+f(c^k v/u)\,I)\,L^{[k]}_1(v)\,(R^{-1}-f(v/u)\,I)}{1-f(v/u)(f(v/u)+q-\qq)}.$$

Now,  we apply the operator $\mathrm{Tr}_{R(2)}$ to this relation. Thus, we get
$$ \frac{(\mathrm{Tr}_R\, L^{[k]}(v))\, L_1(u)}{1-f(u/c^{k-1}v)(f(u/c^{k-1}v)+q-\qq)}-\frac{L_1(u)\, (\mathrm{Tr}_R\, L^{[k]}(v))}{1-f(v/u)(f(v/u)+q-\qq)}=$$
$$\frac{L_1(u)\, L^{[k]}_1(v)(-f(v/u)+f(c^k v/u)(1-\al(q-\qq))-f(v/u)f(c^k v/u)\al)}{1-f(v/u)(f(v/u)+q-\qq)}-$$
\be \frac{L^{[k]}_1(v)\,L_1(u)(-f(u/c^{k-1}v)+f(cu/v)(1-\al(q-\qq))-f(u/c^{k-1}v)f(cu/v)\al)}{1-f(u/c^{k-1}v)(f(u/c^{k-1}v)+q-\qq)}. \label{diff} \ee

Here, we used the equality
$$\mathrm{Tr}_{R(2)} R^{-1}= \mathrm{Tr}_{R(2)} R-(q-\qq) \mathrm{Tr}_{R(2)} I= (1-(q-\qq)\al) I. $$

If $c\not=1$ and $k\geq 2$ then the denominators of the summands in the left hand side of the
relation (\ref{diff}) are distinct. This prevents  us from estimating the difference
\be
(\mathrm{Tr}_R\, L^{[k]}(v))\, L_1(u)-L_1(u)\, (\mathrm{Tr}_R\, L^{[k]}(v)). \label{di}
\ee
Thus, we are not able to argue either for or against the centrality of the  quantum power sums $\mathrm{Tr}_R\, L^{[k]}(v),\, k\geq 2$.

Nevertheless, if we put $c=1$ in the algebra defined by (\ref{RS}) and  (\ref{RRR}), it is not difficult to get a condition on the symmetry $R$ entailing the centrality of the higher
 quantum power sums $\mathrm{Tr}_RL^k(u)$\footnote{Note that if $c=1$ then $L^{[k]}(v)=L^k(v)$, i.e. the quantum powers of $L(u)$ are equal to the usual powers. However,
observe that the trace in the definition of the quantum power sums  is not usual.}. Indeed, in this case the mentioned denominators in  (\ref{diff}) are equal to each other as
a consequence of the explicit form of the function $f(x)$ (\ref{fff}):
$$
f(u/v)(f(u/v)+q-q^{-1}) = f(v/u)(f(v/u)+q-q^{-1})  = (q-q^{-1})^2\,\frac{uv}{(u-v)^2}.
$$

Now, due to this fact we can  estimate the difference (\ref{di}).
After some transformations we get the following equality
\be
(\mathrm{Tr}_R\, L^{k}(v))\, L_1(u)-L_1(u)\, (\mathrm{Tr}_R\, L^{k}(v))=\Big(L_1(u)\, L^{k}(v)-L^{k}(v)\,L_1(u)\Big)\,\frac{\alpha\,(q-q^{-1})uv}{(u-v)^2},
\label{diff-est}
\ee
where $\alpha = q^{n-m}(m-n)_q\not=0$  (recall that $(m|n)$ is the bi-rank of $R$). If $m\not=n$, then $\al\not=0$. Consequently,  the right hand
side of this equality does not vanish. This entails that none of the quantum power sums in the corresponding algebra is central.

However, if $m=n$, i.e. the bi-rank of the Hecke symmetry $R$ is $(m|m)$, the right hand side of (\ref{diff-est}) vanishes.
Note that in this case the critical value of the charge is just  $c=q^{2(m-n)} = 1$. Consequently,  we arrive to the following conclusion.

\begin{proposition}
If the bi-rank of a skew-inventible Hecke symmetry $R$ is $(m|m)$ then in the algebra defined by {\rm (\ref{RS})} with $f(x)$ defined by {\rm (\ref{fff})}
the critical value of the charge equals 1. For this value of the charge all quantum power sums $\mathrm{Tr}_R\, L^{[k]}(v)=\mathrm{Tr}_R\, L^{k}(v)$ are central.
\end{proposition}

Note that there are generalized RE algebras with $R$ of bi-rank $(m|m)$ which are not of RS type with central quantum power sums $\mathrm{Tr}_R\, L^{[k]}(v)=\mathrm{Tr}_R\, L^k(v)$.

In conclusion, we want to observe the following. First, the braided Yangians, introduced by ourselves recently (see \cite{GS, GST}), are particular
cases of the generalized RE algebras: for them $g_2=g_3=0$. (Recall that in the present note we assume all $g_k$ to be nontrivial.)
In these algebras quantum analogs of matrix power and power sums are also well-defined.
However, the latter elements are not central. The only central element in each of them is the quantum determinant $\det_R L(u)$, or more precisely, the coefficients of
its expansion in a Laurent series in $u$.
(Also, recall  that in the present note we do not consider similar expansions of elements of the generalized RE algebras.) Besides,
in the braided Yangians quantum analogs of other symmetric polynomials are well-defined. We know no their reasonable definition in the generalized RE algebra.

Second,  note that the deformation property of the generalized RE algebras (i.e. whether this algebra is a deformation of that $U(\widehat{gl(N)})$) is not clear. According the the paper \cite{DF} the double RTT algebra from  \cite{RS} is isomorph
to  the affine QG $U_q(\widehat{gl(N)})$. However, the map taking the RE algebra of RS type to this QG has a kernel. This prevents one from concluding on the deformation property
of the former algebra.


\begin{thebibliography}{MRZ}

\bibitem[C]{C} I.Cherednik {\em Factorizing particles on a half line, and root systems} (Russian), Teoret. Mat. Fiz. 61 (1984), 35--44.

\bibitem[DF]{DF}  J.Ding, I.Frenkel {\em Isomorphism of two realization of quantum affine algebra $U_q(\widehat{gl(n)})$}, Commun. Math. Phys. 156 (1993), 277--300.

\bibitem[FRT]{FRT} L.Faddeev, N.Reshetikhin, L.Takhtadzhyan {\em Quantization of Lie groups and Lie algebras}, Leningrad Math. J. 1 (1989), 193--225.

\bibitem[G]{G} D.Gurevich {\em Algebraic aspects of Quantum Yang-Baxter equation}, Leningrad Math. Journal 2:4 (1990), 119--148.

\bibitem[GPS]{GPS} D.Gurevich, P.Pyatov, P.Saponov {\em Representation theory of (modified) Reflection Equation Algebra of the $GL(m|n)$ type},
 Algebra and Analysis, {\bf 20} (2008), 70--133.

\bibitem[GS]{GS} D.Gurevich, P.Saponov {\em Braided Yangians}, arXiv:1612.05929.

\bibitem[GST]{GST} D.Gurevich, P.Saponov, D.Talalaev {\em Drinfeld-Sokolov reduction in quantum algebras}, arXiv:1710.01806.

\bibitem[IOP]{IOP} A.Isaev, O.Ogievetsky, P.Pyatov {\em On quantum matrix algebras satisfying the Cayley-Hamilton-Newton identities}, J. Phys. A 32 (1999), no. 9, L115--L121.

\bibitem[KM]{KM} S.Ko\v{z}i\'{c}, A.Molev {\em Center of the quntum affine vertex algebra associated with trigonometrical $R$-matrix}, Journal of Physics A: Mathematical and General 50
(32), 1-21.

\bibitem[KS]{KS} P.Kulish, E.Sklyanin {\em Algebraic structures related to reflection equation}, J.Phys. A 25 (1992)  5963--5975.

\bibitem[O]{O} O.Ogievetsky, {\em Uses of Quantum Spaces}, 3rd cycle. Bariloche (Argentine), 2000, pp.72, cel-00374419.


\bibitem[RS]{RS} N.Reshetikhin, M.Semenov-Tian-Shansky, {\em Central extensions of quantum current groups},  Lett. Math. Phys. 19 (1990), no. 2, pp. 133–142.

\bibitem[S]{S} E.Sklyanin {\em Boundary condition for integrable quantum system}, J.Phys. A 21 (1988), 2375--2389.

\bibitem[T]{T} D.Talalaev {\em Quantum Gaudin system}, Func. Anal. Appl. 40 (2006), no 1, 73-77.
\end{thebibliography}
\end{document}